\documentclass{amsart}
\usepackage{amstext}
\usepackage{amsthm}
\usepackage{amssymb}
\usepackage{esint}

\makeatletter

\numberwithin{equation}{section}
\numberwithin{figure}{section}
\theoremstyle{plain}
\newtheorem{thm}{Theorem}
  \theoremstyle{plain}
  \newtheorem{cor}[thm]{Corollary}

\allowdisplaybreaks
\usepackage[sort&compress,square,numbers]{natbib}

\makeatother

\global\long\def\relphantom#1{\mathrel{\phantom{{#1}}}}

\begin{document}

\title[Some identities of Chebyshev polynomials]{Some identities of Chebyshev polynomials arising from non-linear
differential equations}

\author{Taekyun Kim}
\address{ Department of Mathematics, College of Science, Tianjin Polytechnic University, Tianjin City, 300387, China\\
Department of Mathematics, Kwangwoon University, Seoul 139-701, Republic
of Korea}
\email{tkkim@kw.ac.kr}

\author{Dae San Kim}
\address{Department of Mathematics, Sogang University, Seoul 121-742, Republic
of Korea}
\email{dskim@sogang.ac.kr}

\author{Jong-Jin Seo}
\address{Department of Applied Mathematics, Pukyong National  University, Pusan, Republic of Korea}
\email{seo2011@pknu.ac.kr}

\author{Dmitry V. Dolgy}
\address{School of Natural Sciences, Far Eastern Federal University, Vladivostok, Russia}
\email{dvdolgy@pknu.ac.kr}

\begin{abstract}
In this paper, we investigate some properties of Chebyshev polynomials
arising from non-linear differential equations. From our investigation,
we derive some new and interesting identities on Chebyshev polynomials.
\end{abstract}


\maketitle

\section{Introduction}

As is well known, the Chebyshev polynomials of the first kind, $T_{n}\left(x\right)$,
$\left(n\ge0\right)$, are defined by the generating function
\begin{equation}
\frac{1-t^{2}}{1-2xt+t^{2}}=\sum_{n=0}^{\infty}T_{n}\left(x\right)\frac{t^{n}}{n!},\quad\left(\text{see \cite{key-1,key-3,key-5,key-7,key-17,key-21}}\right).\label{eq:1}
\end{equation}

The higher-order Chebyshev polynomials are given by the generating
function
\begin{equation}
\left(\frac{1-t^{2}}{1-2xt+t^{2}}\right)^{\alpha}=\sum_{n=0}^{\infty}T_{n}^{\left(\alpha\right)}\left(x\right)t^{n},\label{eq:2}
\end{equation}
and Chebyshev polynomials of the second kind are denoted by $U_{n}$
and given by generating function
\begin{equation}
\frac{1}{1-2xt+t^{2}}=\sum_{n=0}^{\infty}U_{n}\left(x\right)t^{n},\quad\left(\text{see \cite{key-1,key-8,key-12,key-17}}\right).\label{eq:3}
\end{equation}

The higher-order Chebyshev polynomials of the second kind are also
defined by
\begin{equation}
\left(\frac{1}{1-2xt+t^{2}}\right)^{\alpha}=\sum_{n=0}^{\infty}U_{n}^{\left(\alpha\right)}\left(x\right)t^{n}.\label{eq:4}
\end{equation}

The Chebyshev polynomials of the third kind are defined by the generating
function
\begin{equation}
\frac{1-t}{1-2xt+t^{2}}=\sum_{n=0}^{\infty}V_{n}\left(x\right)t^{n},\quad\left(\text{see \cite{key-1,key-7,key-8,key-17}}\right).\label{eq:5}
\end{equation}
and the higher-order Chebyshev polynomials of the third kind are also
given by the generating function
\begin{equation}
\left(\frac{1-t}{1-2xt+t^{2}}\right)^{\alpha}=\sum_{n=0}^{\infty}V_{n}^{\left(\alpha\right)}\left(x\right)t^{n}.\label{eq:6}
\end{equation}

Finally, we introduce the Chebyshev polynomials of the fourth kind
defined by the generating function
\begin{equation}
\frac{1+t}{1-2xt+t^{2}}=\sum_{n=0}^{\infty}W_{n}\left(x\right)t^{n}.\label{eq:7}
\end{equation}

The higher-order Chebyshev polynomials of the fourth kind are defined
by
\begin{equation}
\left(\frac{1+t}{1-2xt+t^{2}}\right)^{\alpha}=\sum_{n=0}^{\infty}W_{n}^{\left(\alpha\right)}\left(x\right)t^{n}.\label{eq:8}
\end{equation}

It is well known that the Legendre polynomials are defined by the
generating function
\begin{equation}
\frac{1}{\sqrt{1-2xt+t^{2}}}=\sum_{n=0}^{\infty}p_{n}\left(x\right)t^{n},\quad\left(\text{see \cite{key-2,key-19}}\right).\label{eq:9}
\end{equation}

Chebyshev polynomials are important in approximation theory because
the roots of the Chebyshev polynomials of the first kind, which are
also called Chebyshev nodes, are used as nodes in polynomial nodes
(see \cite{key-18}).

The Chebyshev polynomials of the first kind and of the second kind are
solutions of the following Chebyshev differential equations
\begin{equation}
\left(1-x^{2}\right)y^{\prime\prime}-xy^{\prime}+n^{2}y=0,\label{eq:10}
\end{equation}
and
\begin{equation}
\left(1-x^{2}\right)y^{\prime\prime}-3xy^{\prime}+n\left(n+2\right)y=0.\label{eq:11}
\end{equation}

These equations are special cases of the Strum-Liouville differential
equation (see \cite{key-1,key-2,key-3}).

The Chebyshev polynomials of the first kind can be defined by the
contour integral
\begin{equation}
T_{n}\left(z\right)=\frac{1}{4\pi i}\oint\frac{\left(1-t^{2}\right)}{1-2tz+t^{2}}t^{-n-1}dt,\label{eq:12}
\end{equation}
where the contour encloses the origin and is traversed in a counterclockwise
direction (see \cite{key-1,key-18,key-21}). The formula for $T_{n}\left(x\right)$
is given by
\begin{equation}
T_{n}\left(x\right)=\sum_{m=0}^{\left[\frac{n}{2}\right]}\binom{n}{2m}x^{n-2m}\left(x^{2}-1\right)^{m}.\label{eq:13}
\end{equation}

From (\ref{eq:3}), we note that
\begin{equation}
2\left(x-t\right)\left(1-2xt+t^{2}\right)^{-2}=\sum_{n=0}^{\infty}nU_{n}\left(x\right)t^{n-1}.\label{eq:14}
\end{equation}

Thus, by (\ref{eq:14}), we get
\begin{equation}
\left(2xt-2t^{2}\right)\left(1-2xt+t^{2}\right)^{-2}=\sum_{n=0}^{\infty}nU_{n}\left(x\right)t^{n}.\label{eq:15}
\end{equation}

From (\ref{eq:3}) and (\ref{eq:15}), we can derive the following
equation:
\begin{align}
\frac{\left(2xt-2t^{2}\right)+\left(1-2xt+t^{2}\right)}{\left(1-2xt+t^{2}\right)^{2}} & =\frac{1-t^{2}}{\left(1-2xt+t^{2}\right)^{2}}\label{eq:16}\\
 & =\sum_{n=0}^{\infty}\left(n+1\right)U_{n}\left(x\right)t^{n}.\nonumber
\end{align}

Note that
\begin{align}
 & \frac{1-t^{2}}{\left(1-2xt+t^{2}\right)^{2}}\label{eq:17}\\
 & =\left(\frac{1-t^{2}}{1-2xt+t^{2}}\right)\left(\frac{1}{1-2xt+t^{2}}\right)\nonumber \\
 & =\left(\sum_{l=0}^{\infty}T_{l}\left(x\right)t^{l}\right)\left(\sum_{m=0}^{\infty}U_{m}\left(x\right)t^{m}\right)\nonumber \\
 & =\sum_{n=0}^{\infty}\left(\sum_{l=0}^{n}T_{l}\left(x\right)U_{n-l}\left(x\right)\right)t^{n}.\nonumber
\end{align}

From (\ref{eq:16}) and (\ref{eq:17}), we have
\[
U_{n}\left(x\right)=\frac{1}{n+1}\sum_{l=0}^{n}T_{l}\left(x\right)U_{n-l}\left(x\right).
\]

The Chebyshev polynomials have been studied by many authors in the several
areas (see \cite{key-1,key-2,key-3,key-4,key-5,key-6,key-7,key-8,key-9,key-10,key-11,key-12,key-13,key-14,key-15,key-16,key-17,key-18,key-19,key-20,key-21}).

In \cite{key-10}, Kim-Kim studied non-linear differential equations
arising from Changhee polynomials and numbers related to Chebyshev poynomials.

In this paper, we study non-linear differential equations arising
from Chebyshev polynomials and give some new and explicit formulas
for those polynomials.

\section{Differential equations arising from Chebyshev polynomials and their
applications}

Let
\begin{equation}
F=F\left(t,x\right)=\frac{1}{1-2tx+t^{2}}.\label{eq:18}
\end{equation}

Then, by (\ref{eq:1}), we get
\begin{equation}
F^{\left(1\right)}=\frac{d}{dt}F\left(t,x\right)=2\left(x-t\right)F^{2}.\label{eq:19}
\end{equation}

From (\ref{eq:19}), we note that
\begin{equation}
2F^{2}=\left(x-t\right)^{-1}F^{\left(1\right)}.\label{eq:20}
\end{equation}

By using (\ref{eq:20}) and (\ref{eq:19}), we obtain the following
equations:
\begin{align}
2^{2}\cdot2F^{3} & =\left(x-t\right)^{-3}F^{\left(1\right)}+\left(x-t\right)^{-2}F^{\left(2\right)},\label{eq:20-1}\\
2^{3}\cdot2\cdot3F^{4} & =3\left(x-t\right)^{-5}F^{\left(1\right)}+3\left(x-t\right)^{-4}F^{\left(2\right)}+\left(x-t\right)^{-3}F^{\left(3\right)}\label{eq:21}
\end{align}
and
\begin{align}
2^{4}\cdot2\cdot3\cdot4F^{5} & =3\cdot5\left(x-t\right)^{-6}F^{\left(1\right)}+3\cdot5\left(x-t\right)^{-6}F^{\left(2\right)}\label{eq:22}\\
 & \relphantom =+\left(3\cdot2\right)\left(x-t\right)^{-5}F^{\left(3\right)}+\left(x-t\right)^{-4}F^{\left(4\right)},\nonumber
\end{align}
where
\[
F^{N}=\underset{N-\text{times}}{\underbrace{F\times\cdots\times F}}\quad\text{and}\quad F^{\left(N\right)}=\left(\frac{d}{dt}\right)^{N}F\left(t,x\right).
\]

Continuing this process, we set
\begin{equation}
2^{N}N!F^{N+1}=\sum_{i=1}^{N}a_{i}\left(N\right)\left(x-t\right)^{i-2N}F^{\left(i\right)},\label{eq:23}
\end{equation}
where $N\in\mathbb{N}$.

From (\ref{eq:23}), we note that
\begin{align}
 & \relphantom =2^{N}N!F^{N}\left(N+1\right)F^{\left(1\right)}\label{eq:24}\\
 & =\sum_{i=1}^{N}a_{i}\left(N\right)\left(2N-i\right)\left(x-t\right)^{i-2N-1}F^{\left(i\right)}+\sum_{i=1}^{N}a_{i}\left(N\right)\left(x-t\right)^{i-2N}F^{\left(i+1\right)}.\nonumber
\end{align}

By (\ref{eq:19}) and (\ref{eq:24}), we get
\begin{align}
 & \relphantom =2^{N}N!\left(N+1\right)F^{N}\left(2\left(x-t\right)F^{2}\right)\label{eq:25}\\
 & =\sum_{i=1}^{N}a_{i}\left(N\right)\left(2N-i\right)\left(x-t\right)^{i-2N-1}F^{\left(i\right)}\nonumber \\
 & \relphantom =+\sum_{i=1}^{N}a_{i}\left(N\right)\left(x-t\right)^{i-2N}F^{\left(i+1\right)}.\nonumber
\end{align}

Thus, from (\ref{eq:25}), we have
\begin{align}
 & \relphantom =2^{N+1}\left(N+1\right)!F^{N+2}\label{eq:26}\\
 & =\sum_{i=1}^{N}a_{i}\left(N\right)\left(2N-i\right)\left(x-t\right)^{i-2\left(N+1\right)}F^{\left(i\right)}\nonumber \\
 & \relphantom =+\sum_{i=2}^{N+1}a_{i-1}\left(N\right)\left(x-t\right)^{i-2\left(N+1\right)}F^{\left(i\right)}.\nonumber
\end{align}

On the other hand, by replacing $N$ by $N+1$, in (\ref{eq:23}),
we get
\begin{equation}
2^{N+1}\left(N+1\right)!F^{N+2}=\sum_{i=1}^{N+1}a_{i}\left(N+1\right)\left(x-t\right)^{i-2\left(N+1\right)}F^{\left(i\right)}.\label{eq:27}
\end{equation}

Comparing the coefficients on both sides of (\ref{eq:26}) and
(\ref{eq:27}), we have
\begin{align}
a_{1}\left(N+1\right) & =\left(2N-1\right)a_{1}\left(N\right),\label{eq:28}\\
a_{N+1}\left(N+1\right) & =a_{N}\left(N\right),\label{eq:29}
\end{align}
and
\begin{equation}
a_{i}\left(N+1\right)=a_{i-1}\left(N\right)+\left(2N-i\right)a_{i}\left(N\right),\quad\left(2\le i\le N\right).\label{eq:30}
\end{equation}

Moreover, by (\ref{eq:20-1}) and (\ref{eq:23}), we get
\begin{equation}
2F^{2}=\left(x-t\right)^{-1}F^{\left(1\right)}=a_{1}\left(1\right)\left(x-t\right)^{-1}F^{\left(1\right)}.\label{eq:31}
\end{equation}

By comparing the coefficients on both sides of (\ref{eq:31}),
we get
\begin{equation}
a_{1}\left(1\right)=1.\label{eq:32}
\end{equation}

Now, by (\ref{eq:28}) and (\ref{eq:32}), we have
\begin{align}
a_{1}\left(N+1\right) & =\left(2N-1\right)a_{1}\left(N\right)\label{eq:33}\\
 & =\left(2N-1\right)\left(2N-3\right)a_{1}\left(N-1\right)\nonumber \\
 & =\left(2N-1\right)\left(2N-3\right)\left(2N-5\right)a_{1}\left(N-2\right)\nonumber \\
 & \vdots\nonumber \\
 & =\left(2N-1\right)\left(2N-3\right)\left(2N-5\right)\cdots1\cdot a_{1}\left(1\right)\nonumber \\
 & =\left(2N-1\right)!!,\nonumber
\end{align}
where $\left(2N-1\right)!!$ is Arfken's double factorial.

From (\ref{eq:29}), we easily note that
\begin{equation}
a_{N+1}\left(N+1\right)=a_{N}\left(N\right)=\cdots=a_{1}\left(1\right)=1.\label{eq:34}
\end{equation}

For $2\le i\le N$, from (\ref{eq:30}), we can derive the following
equation:
\begin{align}
a_{i}\left(N+1\right) & =a_{i-1}\left(N\right)+\left(2N-i\right)a_{i}\left(N\right)\label{eq:35}\\
 & =a_{i-1}\left(N\right)+\left(2N-i\right)a_{i-1}\left(N-1\right)+\left(2N-i\right)\left(2N-2-i\right)a_{i}\left(N-1\right)\nonumber \\
 & \vdots\nonumber \\
 & =\sum_{k=0}^{N-i}\left(\prod_{l=0}^{k-1}\left(2\left(N-l\right)-i\right)\right)a_{i-1}\left(N-k\right)+\prod_{l=0}^{N-i}\left(2\left(N-l\right)-i\right)a_{i}\left(i\right)\nonumber \\
 & =\sum_{k=0}^{N-i}2^{k}\left(N-\frac{i}{2}\right)_{k}a_{i-1}\left(N-k\right)+2^{N-i+1}\left(N-\frac{i}{2}\right)_{N-i+1}\nonumber \\
 & =\sum_{k=0}^{N-i+1}2^{k}\left(N-\frac{i}{2}\right)_{k}a_{i-1}\left(N-k\right),\nonumber
\end{align}
where $\left(x\right)_{n}=x\left(x-1\right)\cdots\left(x-n+1\right)$,
$\left(n\ge1\right)$ and $\left(x\right)_{0}=1$.

As the above is also valid for $i=N+1$, by (\ref{eq:35}), we get
\begin{equation}
a_{i}\left(N+1\right)=\sum_{k=0}^{N+1-i}2^{k}\left(N-\frac{i}{2}\right)_{k}a_{i-1}\left(N-k\right),\label{eq:36}
\end{equation}
where $2\le i\le N+1$.

Now, we give an explicit expression for $a_{i}\left(N+1\right)$.

From (\ref{eq:33}) and (\ref{eq:36}), we can derive the following
equations:
\begin{align}
a_{2}\left(N+1\right) & =\sum_{k_{1}=0}^{N-1}2^{k_{1}}\left(N-\frac{2}{2}\right)_{k_{1}}a_{1}\left(N-k_{1}\right)\label{eq:37}\\
 & =\sum_{k_{1}=0}^{N-1}2^{k_{1}}\left(N-\frac{2}{2}\right)_{k_{1}}\left(2\left(N-k_{1}-1\right)-1\right)!!,\nonumber
\end{align}
\begin{align}
a_{3}\left(N+1\right) & =\sum_{k_{2}=0}^{N-2}2^{k_{2}}\left(N-\frac{3}{2}\right)_{k_{2}}a_{2}\left(N-k_{2}\right)\label{eq:38}\\
 & =\sum_{k_{2}=0}^{N-2}\sum_{k_{1}=0}^{N-2-k_{2}}2^{k_{1}+k_{2}}\left(N-\frac{3}{2}\right)_{k_{2}}\left(N-k_{2}-\frac{4}{2}\right)_{k_{1}}\left(2\left(N-2-k_{1}-k_{2}\right)-1\right)!!,\nonumber
\end{align}
and
\begin{align}
a_{4}\left(N+1\right) & =\sum_{k_{3}=0}^{N-3}2^{k_{3}}\left(N-\frac{4}{2}\right)_{k_{3}}a_{3}\left(N-k_{3}\right)\label{eq:39}\\
 & =\sum_{k_{3}=0}^{N-3}\sum_{k_{2}=0}^{N-3-k_{3}}\sum_{k_{1}=0}^{N-3-k_{3}-k_{2}}2^{k_{1}+k_{2}+k_{3}}\left(N-\frac{4}{2}\right)_{k_{3}}\left(N-k_{3}-\frac{5}{2}\right)_{k_{2}}\left(N-k_{3}-k_{2}-\frac{6}{2}\right)_{k_{1}}\nonumber \\
 & \relphantom =\times\left(2\left(N-3-k_{1}-k_{2}-k_{3}\right)-1\right)!!.\nonumber
\end{align}

Thus, we see that, for $2\le i\le N+1$,
\begin{align}
a_{i}\left(N+1\right) & =\sum_{k_{i-1}=0}^{N-i+1}\sum_{k_{i-2}=0}^{N-i+1-k_{i-1}}\cdots\sum_{k_{1}=0}^{N-i+1-k_{i-1}-\cdots-k_{2}}2^{\sum_{j=1}^{i-1}k_{j}}\label{eq:40}\\
 & \relphantom =\times\prod_{j=2}^{i}\left(N-\sum_{l=j}^{i-1}k_{l}-\frac{2i-j}{2}\right)_{k_{j-1}}\left(2\left(N-i+1-\sum_{j=1}^{i-1}k_{j}\right)-1\right)!!.\nonumber
\end{align}

Therefore, we obtain the following theorem.
\begin{thm}
\label{thm:1} The nonlinear differential equations
\[
2^{N}N!F^{N+1}=\sum_{i=1}^{N}a_{i}\left(N\right)\left(x-t\right)^{i-2N}F^{\left(i\right)},\quad\left(N\in\mathbb{N}\right)
\]
has a solution $F=F\left(t,x\right)=\frac{1}{1-2tx+t^{2}}$, where
\begin{align*}
a_{1}\left(N\right) & =\left(2N-3\right)!!,\\
a_{i}\left(N\right) & =\sum_{k_{i-1}=0}^{N-i}\sum_{k_{i-2}=0}^{N-i-k_{i-1}}\cdots\sum_{k_{1}=0}^{N-i-k_{i-1}-\cdots-k_{2}}2^{\sum_{j=1}^{i-1}k_{j}}\\
 & \relphantom =\times\prod_{j=2}^{i}\left(N-\sum_{l=j}^{i-1}k_{l}-\frac{2i+2-j}{2}_{k_{j-1}}\right)\left(2\left(N-i-\sum_{j=1}^{i-1}k_{j}\right)-1\right)!!
\end{align*}
$\left(2\le i\le N\right)$.
\end{thm}
From (\ref{eq:3}) and (\ref{eq:9}), we note that
\begin{align}
 & \sum_{n=0}^{\infty}U_{n}\left(x\right)t^{n}\label{eq:41}\\
 & =\frac{1}{1-2xt+t^{2}}\nonumber \\
 & =\left(\frac{1}{\sqrt{1-2xt+t^{2}}}\right)^{2}\nonumber \\
 & =\left(\sum_{l=0}^{\infty}p_{l}\left(x\right)t^{l}\right)\left(\sum_{m=0}^{\infty}p_{m}\left(x\right)t^{m}\right)\nonumber \\
 & =\sum_{n=0}^{\infty}\left(\sum_{l=0}^{n}p_{l}\left(x\right)p_{n-l}\left(x\right)\right)t^{n}.\nonumber
\end{align}

Thus, from (\ref{eq:41}), we have
\begin{equation}
U_{n}\left(x\right)=\sum_{l=0}^{n}p_{l}\left(x\right)p_{n-l}\left(x\right).\nonumber
\end{equation}

From (\ref{eq:4}), we obtain
\begin{equation}
2^{N}N!F^{N+1}=2^{N}N!\sum_{n=0}^{\infty}U_{n}^{\left(N+1\right)}\left(x\right)t^{n}.\label{eq:42}
\end{equation}

On the other hand, by Theorem \ref{thm:1}, we get
\begin{align}
2^{N}N!F^{N+1} & =\sum_{i=1}^{N}a_{i}\left(N\right)\left(x-t\right)^{i-2N}F^{\left(i\right)}\label{eq:43}\\
 & =\sum_{i=1}^{N}a_{i}\left(N\right)\left(\sum_{m=0}^{\infty}\binom{2N+m-i-1}{m}x^{i-2N-m}t^{m}\right)\left(\sum_{l=0}^{\infty}U_{i+l}\left(x\right)\left(l+i\right)_{i}t^{l}\right)\nonumber \\
 & =\sum_{i=1}^{N}a_{i}\left(N\right)\sum_{n=0}^{\infty}\left\{ \sum_{l=0}^{n}\binom{2N+n-l-i-1}{n-l}x^{i-2N-n+l}U_{l+i}\left(x\right)\left(l+i\right)_{i}\right\} t^{n}\nonumber \\
 & =\sum_{n=0}^{\infty}\left\{ \sum_{i=1}^{N}a_{i}\left(N\right)\sum_{l=0}^{n}\binom{2N+n-l-i-1}{n-l}x^{i+l-2N-n}U_{i+l}\left(x\right)\left(l+i\right)_{i}\right\} t^{n}.\nonumber
\end{align}

Comparing the coefficients on the both sides of (\ref{eq:42}) and (\ref{eq:43}),
we obtain the following theorem.
\begin{thm}
\label{thm:2} For $N\in\mathbb{N}$, and $n\in\mathbb{N}\cup\left\{ 0\right\} $,
the following identity holds.

\[
U_{n}^{\left(N+1\right)}\left(x\right)=\frac{1}{2^{N}N!}\sum_{i=1}^{N}a_{i}\left(N\right)\sum_{l=0}^{n}\binom{2N+n-l-i-1}{n-l}U_{l+i}\left(x\right)x^{i+l-2N-n}\left(l+i\right)_{i}.
\]

\end{thm}
The higher-order Legendre polynomials are given by the generating
function
\begin{equation}
\left(\frac{1}{\sqrt{1-2xt+t^{2}}}\right)^{\alpha}=\sum_{n=0}^{\infty}p_{n}^{\left(\alpha\right)}\left(x\right)t^{n}.\label{eq:44}
\end{equation}

Thus, by \ref{eq:4} and (\ref{eq:43}), we get
\begin{align}
 & \sum_{n=0}^{\infty}U_{n}^{\left(\alpha\right)}\left(x\right)t^{n}\label{eq:45}\\
 & =\left(\frac{1}{1-2xt+t^{2}}\right)^{\alpha}\nonumber \\
 & =\left(\frac{1}{\sqrt{1-2xt+t^{2}}}\right)^{2\alpha}\nonumber \\
 & =\left(\sum_{l=0}^{\infty}p_{l}^{\left(\alpha\right)}\left(x\right)t^{l}\right)\left(\sum_{m=0}^{\infty}p_{m}^{\left(\alpha\right)}\left(x\right)t^{m}\right)\nonumber \\
 & =\sum_{n=0}^{\infty}\left(\sum_{l=0}^{n}p_{l}^{\left(\alpha\right)}\left(x\right)p_{n-l}^{\left(\alpha\right)}\left(x\right)\right)t^{n}.\nonumber
\end{align}

From (\ref{eq:45}), we note that
\begin{equation}
U_{n}^{\left(\alpha\right)}\left(x\right)=\sum_{l=0}^{n}p_{l}^{\left(\alpha\right)}\left(x\right)p_{n-l}^{\left(\alpha\right)}\left(x\right).\label{eq:46}
\end{equation}

Therefore, we obtian the following corollaries.
\begin{cor}
\label{cor:3} For $N\in\mathbb{N}$ and $n\in\mathbb{N}\cup\left\{ 0\right\} $,
we have
\begin{align*}
 & \sum_{l=0}^{n}p_{l}^{\left(N+1\right)}p_{n-l}^{\left(N+1\right)}\left(x\right)\\
 & =\frac{1}{2^{N}N!}\sum_{i=1}^{N}a_{i}\left(N\right)\sum_{l=0}^{n}\binom{2N+n-l-i-1}{n-l}U_{l+i}\left(x\right)\left(l+i\right)_{i}x^{i+l-2N-n}.
\end{align*}

\end{cor}

\begin{cor}
\label{cor:4} For $N\in\mathbb{N}$ and $n\in\mathbb{N}$, we have
\begin{align}
& U_{n}^{\left(N+1\right)}\left(x\right)\nonumber\\
& =\frac{1}{2^{N}N!}\sum_{i=1}^{N}a_{i}\left(N\right)\sum_{l=0}^{n}\sum_{j=0}^{l+i}\binom{2N+n-l-i-1}{n-l}x^{i+l-2N-n}\left(l+i\right)_{i}\left(x\right)p_{l+i-j}\left(x\right).\nonumber
\end{align}

\end{cor}
By (\ref{eq:6}), we get
\begin{align}
 & 2^{N}N!F^{N+1}\label{eq:47}\\
 & =2^{N}N!\left(1-t\right)^{-N-1}\left(\frac{1-t}{1-2xt+t^{2}}\right)^{N+1}\nonumber \\
 & =2^{N}N!\left(\sum_{m=0}^{\infty}\binom{N+m}{m}t^{m}\right)\left(\sum_{l=0}^{\infty}V_{l}^{\left(N+1\right)}\left(x\right)t^{l}\right)\nonumber \\
 & =2^{N}N!\sum_{n=0}^{\infty}\left(\sum_{l=0}^{n}\binom{N+n-l}{n-l}V_{l}^{\left(N+1\right)}\left(x\right)\right)t^{n}.\nonumber
\end{align}

On the other hand, by Theorem \ref{thm:1}, we have
\begin{align}
2^{N}N!F^{N+1} & =\sum_{i=1}^{N}a_{i}\left(N\right)\left(x-t\right)^{i-2N}F^{\left(i\right)}\label{eq:48}\\
 & =\sum_{i=1}^{N}a_{i}\left(N\right)\left(x-t\right)^{i-2N}\left(\frac{d}{dt}\right)^{i}\left(\frac{1}{1-t}\cdot\frac{1-t}{1-xt+t^{2}}\right).\nonumber
\end{align}

From Leibniz formula, we note that
\begin{align}
 & \left(\frac{d}{dt}\right)^{i}\left(\frac{1-t}{1-2xt+t^{2}}\cdot\frac{1}{1-t}\right)\label{eq:49}\\
 & =\sum_{l=0}^{i}\binom{i}{l}\left(\left(\frac{d}{dt}\right)^{i-l}\frac{1}{1-t}\right)\left(\left(\frac{d}{dt}\right)^{l}\frac{1-t}{1-2xt+t^{2}}\right)\nonumber \\
 & =\sum_{l=0}^{i}\binom{i}{l}\left(i-l\right)!\left(1-t\right)^{-i+l-1}\left(\frac{d}{dt}\right)^{l}\left(\frac{1-t}{1-2xt+t^{2}}\right)\nonumber \\
 & =\sum_{l=0}^{i}\binom{i}{l}\left(i-l\right)!\sum_{s=0}^{\infty}\binom{i-l+s}{s}t^{s}\sum_{p=0}^{\infty}V_{p+l}\left(x\right)\left(p+l\right)_{l}t^{p}\nonumber \\
 & =\sum_{l=0}^{i}\frac{i!}{l!}\sum_{s=0}^{\infty}\binom{i-l+s}{s}t^{s}\sum_{p=0}^{\infty}V_{p+l}\left(x\right)\left(p+l\right)_{l}t^{p}.\nonumber
\end{align}

By (\ref{eq:48}) and (\ref{eq:49}), we get
\begin{align}
 & 2^{N}N!F^{N+1}\label{eq:50}\\
 & =\sum_{n=0}^{\infty}\left\{ \sum_{i=1}^{N}\sum_{l=0}^{i}a_{i}\left(N\right)\frac{i!}{l!}\sum_{m+s+p=n}\binom{2N+m-i-1}{m}\binom{i-l+s}{s}\right.\nonumber \\
 & \relphantom =\left.\times\left(p+l\right)_{l}x^{i-2N-m}V_{p+l}\left(x\right)\right\} t^{n}.\nonumber
\end{align}

Therefore, by (\ref{eq:47}) and (\ref{eq:50}), we obtain the following
theorem.
\begin{thm}
\label{thm:5} For $N\in\mathbb{N}$ and $n\in\mathbb{N}\cup\left\{ 0\right\} $,
we have the following identity:
\begin{align*}
 & \sum_{l=0}^{n}\binom{N+n-l}{n-l}V_{l}^{\left(N+1\right)}\left(x\right)\\
 & =\frac{1}{2^{N}N!}\sum_{i=1}^{N}\sum_{l=0}^{i}a_{i}\left(N\right)\frac{i!}{l!}\sum_{m+s+p=n}\binom{2N+m-i-1}{m}\binom{i-l+s}{s}\left(p+l\right)_{l}\\
 & \relphantom =\times x^{i-2N-m}V_{p+l}\left(x\right).
\end{align*}

\end{thm}
From (\ref{eq:8}), we note that
\begin{align}
 & 2^{N}N!F^{N+1}\label{eq:51}\\
 & =2^{N}N!\left(1+t\right)^{-N-1}\left(\frac{1+t}{1-2xt+t^{2}}\right)^{N+1}\nonumber \\
 & =2^{N}N!\left(\sum_{m=0}^{\infty}\binom{N+m}{m}\left(-1\right)^{m}t^{m}\right)\left(\sum_{l=0}^{\infty}W_{l}^{\left(N+1\right)}\left(x\right)t^{l}\right)\nonumber \\
 & =2^{N}N!\sum_{n=0}^{\infty}\left(\sum_{l=0}^{n}\left(-1\right)^{n-l}\binom{N+n-l}{n-l}W_{l}^{\left(N+1\right)}\left(x\right)\right)t^{n}.\nonumber
\end{align}

On the other hand, by Theorem \ref{thm:1}, we get
\begin{equation}
2^{N}N!F^{N+1}=\sum_{i=1}^{N}a_{i}\left(N\right)\left(x-t\right)^{i-2N}\left(\frac{d}{dt}\right)^{i}\left\{ \frac{1}{1+t}\cdot\frac{1+t}{1-2xt+t^{2}}\right\} .\label{eq:52}
\end{equation}

Now, we observe that
\begin{align}
 & \left(\frac{d}{dt}\right)^{i}\left\{ \left(\frac{1}{1+t}\right)\left(\frac{1+t}{1-2xt+t^{2}}\right)\right\} \label{eq:53}\\
 & =\sum_{l=0}^{i}\binom{i}{l}\left(-1\right)^{i-l}\left(i-l\right)!\left(\frac{1}{1+t}\right)^{i-l+1}\left(\frac{d}{dt}\right)^{l}\left(\frac{1+t}{1-2xt+t^{2}}\right)\nonumber \\
 & =\sum_{l=0}^{i}\binom{i}{l}\left(-1\right)^{i-l}\left(i-l\right)!\sum_{s=0}^{\infty}\binom{i-l+s}{s}\left(-1\right)^{s}t^{s}\sum_{p=0}^{\infty}W_{p+l}\left(x\right)\left(p+l\right)_{l}t^{p}.\nonumber
\end{align}

From (\ref{eq:52}) and (\ref{eq:53}), we have
\begin{align}
 & 2^{N}N!F^{N+1}\label{eq:54}\\
 & =\sum_{n=0}^{\infty}\left\{ \sum_{i=1}^{N}a_{i}\left(N\right)\sum_{l=0}^{i}\left(-1\right)^{i-l}\frac{i!}{l!}\sum_{m+s+p=n}\left(-1\right)^{s}\binom{2N+m-i-1}{m}\right.\nonumber \\
 & \relphantom =\left.\times\binom{i-l+s}{s}\left(p+l\right)_{l}x^{i-2N-m}W_{p+l}\left(x\right)\right\} t^{n}.\nonumber
\end{align}

Therefore, by (\ref{eq:51}) and (\ref{eq:54}), we obtain the following
theorem.
\begin{thm}
\label{thm:6} For $N\in\mathbb{N}$ and $n\in\mathbb{N}\cup\left\{ 0\right\} $,
the following identity is valid:
\begin{align*}
 & \sum_{l=0}^{n}\left(-1\right)^{n-l}\binom{N+n-l}{n-l}W_{l}^{\left(N+1\right)}\left(x\right)\\
 & =\frac{1}{2^{N}N!}\sum_{i=1}^{N}\sum_{l=0}^{i}\left(-1\right)^{i-l}a_{i}\left(N\right)\frac{i!}{l!}\sum_{m+s+p=n}\left(-1\right)^{s}\binom{2N+m-i-1}{m}\\
 & \relphantom =\times\binom{i-l+s}{s}\left(p+l\right)_{l}x^{i-2N-m}W_{p+l}\left(x\right).
\end{align*}

\end{thm}
From (\ref{eq:1}), we have
\begin{align}
 & 2^{N}N!F^{N+1}\label{eq:55}\\
 & =2^{N}N!\left(\frac{1}{1-t^{2}}\cdot\frac{1-t^{2}}{1-2xt+t^{2}}\right)^{N+1}\nonumber \\
 & =2^{N}N!\left(\frac{1}{1-t}\right)^{N+1}\left(\frac{1}{1+t}\right)^{N+1}\left(\frac{1-t^{2}}{1-2xt+t^{2}}\right)^{N+1}\nonumber \\
 & =2^{N}N!\left(\sum_{l=0}^{\infty}\binom{N+l}{l}t^{l}\right)\left(\sum_{m=0}^{\infty}\binom{m+N}{m}\left(-1\right)^{m}t^{m}\right)\left(\sum_{p=0}^{\infty}T_{p}^{\left(N+1\right)}\left(x\right)t^{p}\right)\nonumber \\
 & =2^{N}N!\sum_{n=0}^{\infty}\left(\sum_{l+m+p=n}\binom{N+l}{l}\binom{m+N}{m}\left(-1\right)^{m}T_{p}^{\left(N+1\right)}\left(x\right)\right)t^{n}.\nonumber
\end{align}

On the other hand, by Theorem \ref{thm:1}, we get
\begin{align}
 & 2^{N}N!F^{N+1}\label{eq:56}\\
 & =\sum_{i=1}^{N}a_{i}\left(N\right)\left(x-t\right)^{i-2N}F^{\left(i\right)}\nonumber \\
 & =\frac{1}{2}\sum_{i=1}^{N}a_{i}\left(N\right)\left(x-t\right)^{i-2N}\left(\frac{d}{dt}\right)^{i}\left\{ \left(\frac{1}{1-t}+\frac{1}{1+t}\right)\frac{1-t^{2}}{1-2xt+t^{2}}\right\} .\nonumber
\end{align}

From Leibniz formula, we note that the following equations:
\begin{align}
 & \left(\frac{d}{dt}\right)^{i}\left\{ \left(\frac{1}{1-t}\right)\cdot\left(\frac{1-t^{2}}{1-2xt+t^{2}}\right)\right\} \label{eq:57}\\
 & =\sum_{l=0}^{i}\binom{i}{l}\left(i-l\right)!\sum_{s=0}^{\infty}\binom{i+s-l}{s}t^{s}\sum_{p=0}^{\infty}T_{p+l}\left(x\right)\left(p+l\right)_{l}t^{p},\nonumber
\end{align}
and
\begin{align}
 & \left(\frac{d}{dt}\right)^{i}\left\{ \left(\frac{1}{1+t}\right)\left(\frac{1-t^{2}}{1-2xt+t^{2}}\right)\right\} \label{eq:58}\\
 & =\sum_{l=0}^{i}\binom{i}{l}\left(i-l\right)!\left(-1\right)^{i-l}\sum_{s=0}^{\infty}\binom{i-l+s}{s}\left(-1\right)^{s}t^{s}\sum_{p=0}^{\infty}T_{p+l}\left(x\right)\left(p+l\right)_{l}t^{p}.\nonumber
\end{align}

By (\ref{eq:56}), (\ref{eq:57}), and (\ref{eq:58}), we obtain
\begin{align}
 & 2^{N}N!F^{N+1}\label{eq:59}\\
 & =\frac{1}{2}\sum_{i=1}^{N}a_{i}\left(N\right)\left(x-t\right)^{i-2N}\sum_{l=0}^{i}\binom{i}{l}\left(i-l\right)!\sum_{s=0}^{\infty}\binom{i+s-l}{s}t^{s}\sum_{k=0}^{\infty}T_{p+l}\left(x\right)\left(p+l\right)_{l}t^{p}\nonumber \\
 & \relphantom =+\frac{1}{2}\sum_{i=1}^{N}a_{i}\left(N\right)\left(x-t\right)^{i-2N}\sum_{l=0}^{i}\binom{i}{l}\left(i-l\right)!\left(-1\right)^{i-l}\sum_{s=0}^{\infty}\binom{i-l+s}{s}\left(-1\right)^{s}t^{s}\nonumber \\
 & \relphantom =\times\sum_{p=0}^{\infty}T_{p+l}\left(x\right)\left(p+l\right)_{l}t^{p}\nonumber \\
 & =\frac{1}{2}\sum_{n=0}^{\infty}\sum_{i=1}^{N}\sum_{l=0}^{i}a_{i}\left(N\right)\frac{i!}{l!}\sum_{m+s+p=n}\binom{2N+m-i-1}{m}\binom{i+s-l}{s}\left(p+l\right)_{l}\nonumber\\
 & \relphantom =\times x^{i-2N-m}T_{p+l}\left(x\right)t^{n}+\frac{1}{2}\sum_{n=0}^{\infty}\sum_{i=1}^{N}\sum_{l=0}^{i}a_{i}\left(N\right)\frac{i!}{l!}\left(-1\right)^{i-l}\nonumber\\
 & \relphantom =\times\sum_{m+s+p=n}\left(-1\right)^{s}\binom{2N+m-i-1}{m}\binom{i+s-l}{s}\left(p+l\right)_{l}x^{i-2N-m}T_{p+l}\left(x\right)t^{n}.\nonumber
\end{align}

Therefore, by (\ref{eq:55}) and (\ref{eq:59}), we obtain the following
theorem.
\begin{thm}
\label{thm:7} For $n\in\mathbb{N}\cup\left\{ 0\right\} $ and $N\in\mathbb{N}$,
we have the following identity
\begin{align*}
 & 2^{N+1}N!\sum_{s+m+p=n}\binom{N+s}{s}\binom{m+N}{m}\left(-1\right)^{m}T_{p}^{\left(N+1\right)}\left(x\right)\\
 & =\sum_{i=1}^{N}\sum_{l=0}^{i}a_{i}\left(N\right)\frac{i!}{l!}\sum_{m+s+p=n}\binom{2N+m-i-1}{m}\binom{i+s-l}{s}\left(p+l\right)_{l}\nonumber\\
 & \relphantom = \times x^{i-2N-m}T_{p+l}\left(x\right)+\sum_{i=1}^{N}\sum_{l=0}^{i}a_{i}\left(N\right)\frac{i!}{l!}\left(-1\right)^{i-l}\sum_{m+s+p=n}\left(-1\right)^{s}\binom{2N+m-i-1}{m}\nonumber\\
 & \relphantom =\times\binom{i+s-l}{s}\left(p+l\right)_{l}x^{i-2N-m}T_{p+l}\left(x\right).
\end{align*}
\end{thm}

Acknowledgements. This paper is supported by grant NO 14-11-00022 of Russian Scientific Fund.

\bibliographystyle{amsplain}
\providecommand{\bysame}{\leavevmode\hbox to3em{\hrulefill}\thinspace}
\providecommand{\MR}{\relax\ifhmode\unskip\space\fi MR }
\providecommand{\MRhref}[2]{%
  \href{http://www.ams.org/mathscinet-getitem?mr=#1}{#2}
}
\providecommand{\href}[2]{#2}

\end{document}